\documentclass{article}
\usepackage{epsfig}
\usepackage{amsmath}
\usepackage{amssymb,latexsym}
\usepackage{a4wide}
\usepackage{makeidx}

\newtheorem{theorem}{Theorem}

\newtheorem{definition}{Definition}
\newtheorem{remark}{Remark}

\def\halmos{\rule{6pt}{6pt}}

\def\det{{\rm det}\, }

\def\vol{{\rm vol}\, }

\def\arccos{{\rm arccos}\, }
\def\arctan{{\rm arctan}\, }

\title{Two-sided bounds for dihedral angle sums of path and $4$-ball tetrahedra}
\author{{\it Sergey Korotov} \\
\centerline{\small Division of Mathematics and Physics, UKK, M\"alardalen University,
Box 883, 721\,23 V\"aster{\aa}s, Sweden}\\
{\small e-mail: {\tt sergey.korotov@mdu.se}}\\ \\
{\it Michal K\v r\'{\i}\v zek} \\
{\small Institute of Mathematics, Czech Academy of Sciences, \v Zitn\'a 25,
115 67 Prague 1, Czech Republic}\\ 
{\small e-mail: {\tt krizek@math.cas.cz}}
}

\begin{document}
\date{\today}
\maketitle

{\bf Abstract:}  A tetrahedron is called a path tetrahedron, if it has three
mutually orthogonal edges that do not intersect at a single point. A tetrahedron is called a $4$-ball tetrahedron, if
there exists
a sphere tangent to all its edges. We derive two-sided tight bounds for
dihedral angle sums of such tetrahedra. In particular, we prove that this sum
lies in the interval $(2\pi, 2.5\pi)$ for path tetrahedra and in $[6\,\arccos\tfrac13, 3\pi)$ for $4$-ball tetrahedra.
Also some of their useful properties are presented.

\smallskip


\smallskip


\section{Introduction}

It was shown by \cite{BraCihKri-2015, BOOK, Gad, Juzuk, Katsuura} that the sum of six dihedral angles between faces
of an arbitrary tetrahedron is not
constant. It varies between $2 \pi$ and $3 \pi$
and these two-sided bounds are tight (optimal), i.e. they cannot be improved, see Theorem~\ref{th:estimate-Gaddum} below stated in 1941.
Since that time a few more results on the dihedral angle sums for various types of tetrahedra and some other polyhedra have been obtained, see
\cite{Hajja-Hayajneh, KorLunVat, KorVat}.

In this paper, we mainly concentrate on a special kind of tetrahedra, the so-called $4$-ball tetrahedra, for which there exists a sphere which touches all six
edges. In Section~\ref{sec:estimation-path-tet}, we derive a simple two-sided bound of the dihedral angle
sum for path tetrahedra (see Theorem~\ref{th:sum-for-path-tets}). These direct generalizations of right triangles
into three-dimensional space are elementary building blocks of all polyhedra, see \cite{KorKri-Dissection}.
Note that tetrahedra with nonobtuse dihedral angles
(cf.~(\ref{eq:two-sided-estimate-path-tet}) and (\ref{eq:two-sided-estim-cube-corner})) have many useful applications in numerical mathematics and
finite element analysis (see \cite[p.~67]{BOOK}).
In Section~\ref{sec:basic-properties-karkas}, we present several basic properties of $4$-ball tetrahedra.
We illustrate that $4$-ball tetrahedra can be very badly shaped, in
particular, all their faces can be obtuse triangles and at the same time three of their dihedral
angles can also be obtuse. Moreover, we address the question of how to divide face-to-face a given
$4$-ball tetrahedron into path tetrahedra so that there are no obtuse angles.
In Section~\ref{sec:estimation-karkas}, we show that the dihedral angle sum
of any $4$-ball tetrahedron is not less than $6\,\arccos\tfrac13\approx 2.35\pi$ 
and is less than $3\pi$, and that each value in this interval is reachable (see Theorem~\ref{th:main-estimate-karkas-tet}).

Throughout the text, all segments and their lengths are denoted for simplicity by the same symbols $a,b,l_1,\dots$

\begin{definition}
  \hspace{-0.15cm}{\bf .}
  \rm
  The {\it tetrahedron}, denoted by the symbol $T$ in what follows, is the convex hull of four points in ${\bf R}^3$
  not lying in the same plane.  The six interior angles between its faces are called {\it dihedral angles}, i.e., each dihedral
  angle is associated with one of the six edges of $T$.
\end{definition}

\begin{definition}
  \hspace{-0.15cm}{\bf .}
  \rm
The {\it dihedral angle sum of the tetrahedron} $T$, denoted by $\Sigma_T$, is the sum of all its six dihedral angles.
\end{definition}

The following two-sided bound for a general tetrahedron was proved by Jarden (Juzuk) \cite{Juzuk} in 1941 (see also \cite{Gad} and \cite{Katsuura}).

\begin{theorem}
    \hspace{-0.15cm}{\bf .}
  For the dihedral angle sum $\Sigma_T$ of a tetrahedron $T$ we have the following estimates
  \begin{equation}
2 \pi < \Sigma_T < 3 \pi,
\label{eq:two-sided-estimate-general-tet}
  \end{equation}
where  both lower and upper bounds are tight.
\label{th:estimate-Gaddum}
\end{theorem}

\section{Two-sided estimation of the dihedral angle sum for path tetrahedra}
\label{sec:estimation-path-tet}

First we present the well-known definition, which generalizes the concept of a right triangle
to three-dimensional space.

\begin{figure}[htbp]
  \centerline{
    \quad
    \psfig{figure=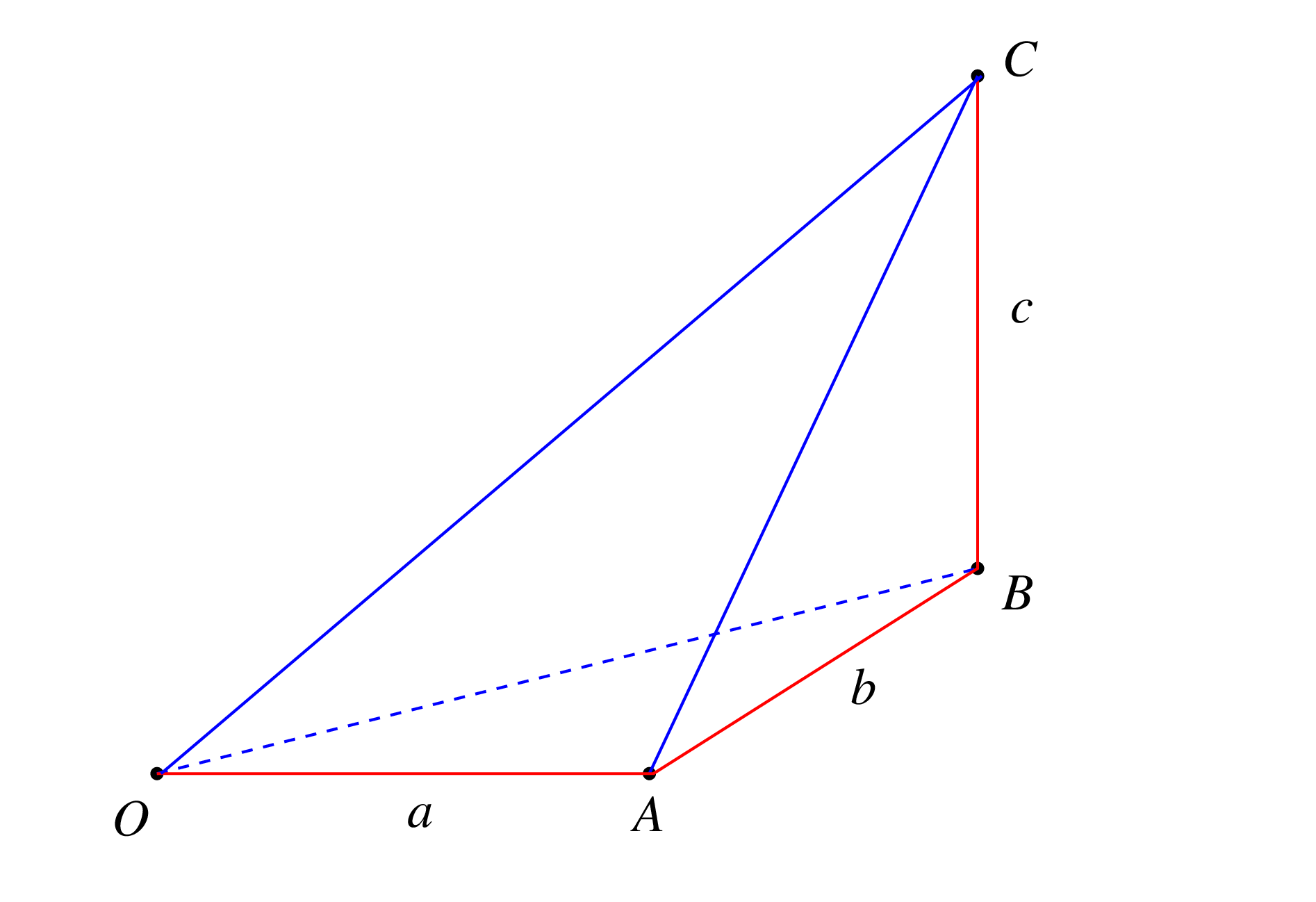, width=8.0cm,height=5.5cm}
  \
    \psfig{figure=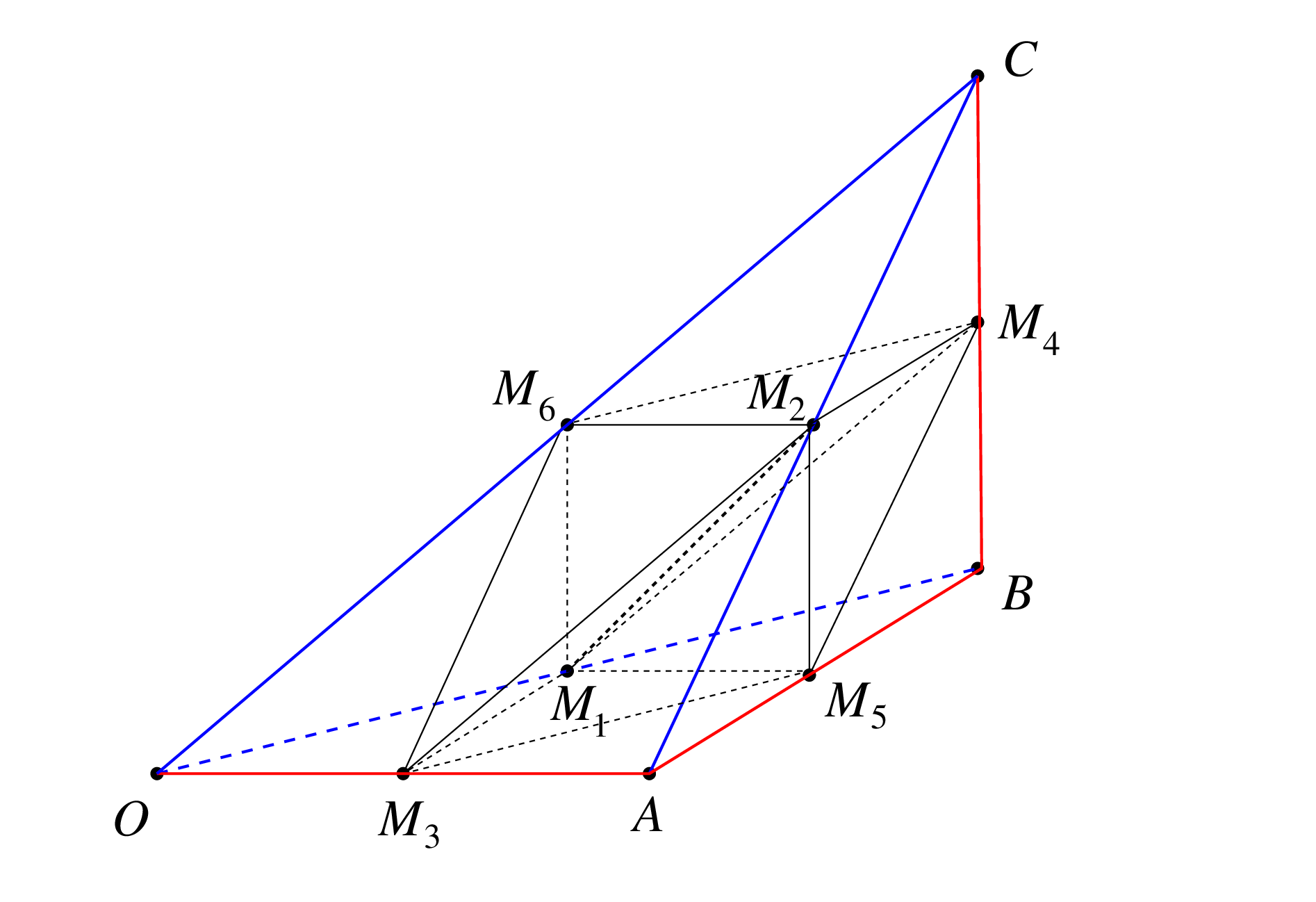, width=8.0cm,height=5.5cm}
}
  \caption{Left: Path tetrahedron $OABC$ with orthogonal edges $a$, $b$, and $c$ forming a path
    (in the sense of graph theory). Its four faces are right triangles and its volume is $abc/6$.
    Right: Using midlines of its faces and the inner diagonal connecting the midpoints $M_1$ and $M_2$
    of the edges $OB$ and $AC$,  it can be partitioned face-to-face into $8$
  path subtetrahedra, see \cite[p.~90]{BOOK}.} 
\label{fig:path-tet}
\end{figure}

\begin{definition}
    \hspace{-0.15cm}{\bf .}
  A tetrahedron $T$ is called a {\it path tetrahedron}, if it has three mutually orthogonal edges that do not intersect at a single point,
  see Figure~\ref{fig:path-tet}.
\end{definition}

Note that all dihedral angles of a path tetrahedron are nonobtuse and all its faces are right triangles, see e.g. \cite[p.~61--62]{BOOK}.

\begin{theorem}
    \hspace{-0.15cm}{\bf .}
Let $T$ be a path tetrahedron. Then
  \begin{equation}
2 \pi < \Sigma_T < 2.5 \pi
\label{eq:two-sided-estimate-path-tet}
  \end{equation}
and each value between these two bounds is attained by some path tetrahedron.
\label{th:sum-for-path-tets}
\end{theorem}

\noindent{\bf P r o o f :} 
The lower and upper bounds hold for any tetrahedron all of whose dihedral
angles are nonobtuse, see \cite[p.~400]{BraCihKri-2015}. So we only demonstrate that each value  in  $(2\pi,2.5\pi)$
is reachable by $\Sigma_T$ of some path simplex $T$. Any path tetrahedron is
completely described by the lengths of three mutually orthogonal edges forming a path. Let $a$,
$b$, $c$ be any
positive numbers and consider a path tetrahedron $T$ with vertices
$O = (0, 0, 0)$, $A = (a, 0, 0)$, $B = (a, b, 0)$, and $C = (a, b, c)$,
see Figure~\ref{fig:path-tet}. The three dihedral angles at edges $OB$, $BA$, and $AC$ are obviously right.
Notice that these edges also form a path.
Without loss of generality we can assume that
$$
b = 1.
$$
Then the angle at $AO$ is $\arctan c$ and the angle at $BC$ is $\arctan
a$. The cosine of the
last dihedral angle at $CO$ is equal to the scalar product of outward
unit normals
with the opposite sign (see e.g. (2.9) of [2]). The outward unit normal to
the face $ACO$ is
$(0, -c , 1)/\sqrt{c^2+1}$ and the outward unit normal to the face $BCO$
is
$(-1,a, 0)/\sqrt{a^2+1}$. Moreover, let $c=a$. Then
$$
\Sigma_T=\Sigma_T(a)=3\frac{\pi}{2} + 2 \, \arctan a +
\arccos\frac{a^2}{a^2+1}
$$
and
$$
\lim_{a \to 0}   \Sigma_T(a)=\frac32\pi+ 0 + \frac12\pi = 2\pi , \qquad
\lim_{a\to\infty}\Sigma_T(a)=\frac32\pi + 2\frac12\pi +0 = 2.5\pi.
$$
Since the function $a\mapsto \Sigma_T(a)$ is continuous (i.e.~having the
Darboux property), the proposed two-sided estimation (\ref{th:sum-for-path-tets})
follows.  \halmos

\section{Basic properties of $4$-ball tetrahedra}
\label{sec:basic-properties-karkas}

\begin{definition}
      \hspace{-0.15cm}{\bf .}
\rm
 A tetrahedron $T$ is called a {\it $4$-ball tetrahedron}, if
 there exists a sphere tangent to all edges of $T$ that will be called a {\it mid-sphere}
 in what follows (see \cite{Matizen-Kvant}). 
\end{definition}

\begin{remark}
      \hspace{-0.15cm}{\bf .}
\rm
In the literature, several other terms for this kind of tetrahedra can
be found -- Crelle's tetrahedron, frame tetrahedron, mid-sphere
tetrahedron, circumscriptible tetrahedron, karkas tetrahedron, quasicircumscribed tetrahedron,
balloon tetrahedron, see e.g. \cite{Court, Crelle, Hajja-2006, Skopets-Ponarin} and the references therein.
\end{remark}

\begin{figure}[htbp]
\centerline{
  \psfig{figure=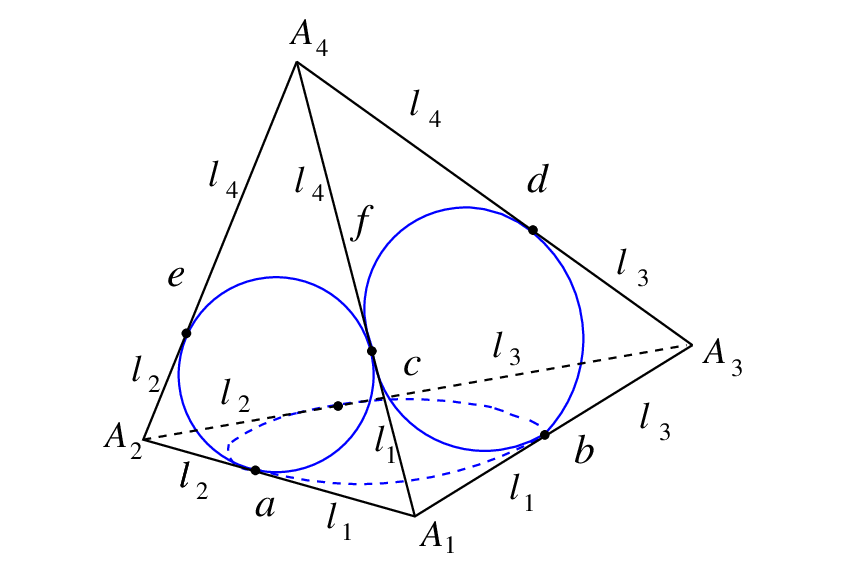, width=8.3cm,height=6cm}
}
\caption{A sketch of a $4$-ball tetrahedron $A_1A_2A_3A_4$. Here $l_i$ are lengths
of tangents from its vertices to the corresponding mid-sphere. The
six touching points are marked by bullets.}
\label{fig:karkas-tangents}
\end{figure}

\begin{remark}
      \hspace{-0.15cm}{\bf .}
      \rm
Consider an arbitrary 4-ball tetrahedron $T$ with vertices
$A_i$, $i\in\{1,2,3,4\}$. According to Theorem~\ref{th:karkas-equivalent-definitions} below, there exist 4
mutually kissing closed balls with centers $A_i$ and radii $\ell_i$.
Then by \cite[p.\,136]{Ponarin}, the radius $\rho$ of the mid-sphere can be computed as follows 
$$
\rho=\frac{2 l_1 l_2 l_3 l_4}{3\,{\rm vol}_3 T},
$$
where ${\rm vol}_3 T$ denotes the volume of $T$. Obviously, $l_i$ is the
length of a tangent from
the vertex $A_i$ of $T$ to its mid-sphere for each $i\in\{1,2,3,4\}$.
Moreover, we can see from Figure~\ref{fig:karkas-tangents}
that the sum of the lengths of any two
opposite edges of any $4$-ball tetrahedron 
is equal to
$$
l_1+l_2+l_3+l_4 = {\rm const}.
$$
Compare this with 1) in Theorem~\ref{th:karkas-equivalent-definitions} below.
\label{rm:karkas-radius}
\end{remark}

The next theorem is very simple, but we will use it several times in
what follows.

\begin{theorem}
        \hspace{-0.15cm}{\bf .}
Let $l_1, l_2, l_3$ be arbitrary positive numbers. Then there exists a triangle with sides
\begin{equation}
  a=l_1+l_2,\qquad b=l_1+l_3,\qquad c=l_2+l_3.
  \label{eq:three}
\end{equation}
Conversely, for an arbitrary triangle with sides $a$, $b$, $c$ there exist uniquely determined positive numbers $l_1$, $l_2$, $l_3$ such that
(\ref{eq:three}) holds, i.e. there exists a one-to-one mapping between the triples $(a, b, c)$ and $(l_1, l_2, l_3)$.
\label{th:auxiliary-theorem}
\end{theorem}

\noindent{\bf P r o o f :}
It is easy to check that all three inequalities required for the existence of a triangle are valid, e.g.,
$$
a+b=l_1+l_2 + l_1+l_3 > l_2+l_3 = c.  
$$

Conversely, for a given triangle with sides $a,b,c$ we put

$$
l_1=\tfrac12(a+b-c), \qquad l_2=\tfrac12(a+c-b), \qquad  l_3=\tfrac12(b+c-a)
$$
to find that (\ref{eq:three}) holds.
The numbers $l_i$ are obviously positive due to triangle inequalities.

Formula (\ref{eq:three})  can be represented by the linear transformation
$$
\left(\begin{matrix}  a \\
                      b \\ 
                      c 
\end{matrix}\right)=
\left(\begin{matrix}  1 & 1 & 0 \\
                      1 & 0 & 1 \\ 
                      0 & 1 & 1 
\end{matrix}\right)
\left(\begin{matrix}  l_1 \\
                      l_2 \\ 
                      l_3 
\end{matrix}\right)
$$
whose inverse transformation reads
$$
\left(\begin{matrix}  l_1 \\
                      l_2 \\ 
                      l_3 
\end{matrix}\right)=\frac12
\left(\begin{matrix}  ~~1 & ~~1 &  -1 \\
                      ~~1 &  -1 & ~~1 \\ 
                       -1 & ~~1 & ~~1 
\end{matrix}\right)
\left(\begin{matrix}  a \\
                      b \\ 
                      c 
\end{matrix}\right).
$$
This implies the required uniqueness.
\halmos

\medskip

In Remark~\ref{rem:no-analog} we demonstrate that an analogous theorem for $4$-ball tetrahedra does not hold. In particular,
there is no unique correspondence between arbitrary positive numbers
$l_i, \, i=1, \dots, 4$, and the lengths of edges of $T$.

\medskip

The next theorem gives several equivalent definitions of a $4$-ball tetrahedron, see \cite{Katsuura-2021} and  \cite[p.\,136]{Ponarin}  
for the corresponding proofs.

\begin{figure}[htbp]
\centerline{
  \psfig{figure=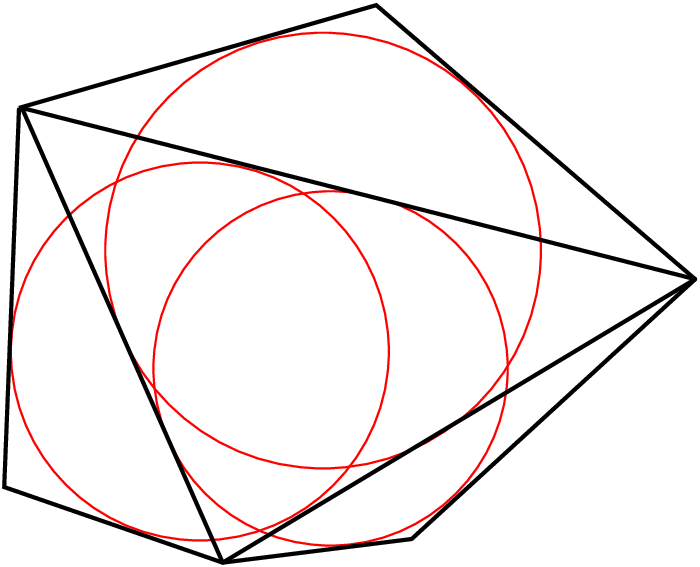, width=6.5cm,height=5.4cm}
  \qquad
    \psfig{figure=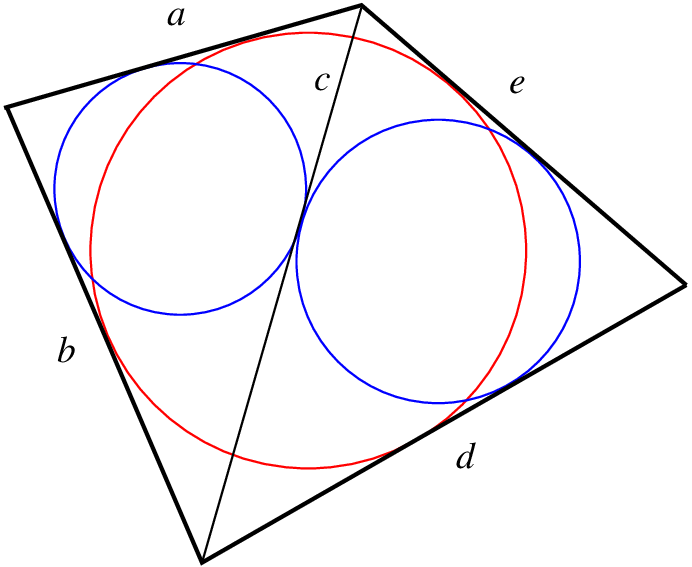, width=6.5cm,height=5.4cm}
}
\caption{Left: Inscribed circles to unfolded quadrilaterals constructed from triangular faces. Right: If there exists an inscribed circle
  to a convex unfolded quadrilateral, then
$a+d = b+e$ (cf.~Remark~\ref{rem:no-analog}). Two smaller circles correspond to the mid-sphere.}
\label{fig:Kvant-Fig-7-8}
\end{figure}

\begin{theorem}
          \hspace{-0.15cm}{\bf .}
A tetrahedron is $4$-ball if and only if any of the following
conditions holds:

\smallskip

1) The sums of lengths of opposite edges are equal.

\smallskip

2) The sums of dihedral angles adjacent to  opposite edges are equal.

\smallskip

3) The circles inscribed in the adjacent faces of the tetrahedron are touching each other in one point
(see Figures~\ref{fig:karkas-tangents} and \ref{fig:Kvant-Fig-7-8} (right)).

\smallskip

4) There exist four mutually externaly tangent balls with
centers at the four vertices.

\smallskip

5) Perpendiculars to the four faces of a tetrahedron emanating from the centers of the circles inscribed in the faces 
intersect at a single point (the center of the mid-sphere).
\label{th:karkas-equivalent-definitions}
\end{theorem}

Consider a $4$-ball tetrahedron and its mid-sphere. Then the three straight lines connecting the touching points of opposite edges intersect at one
point, for the proof see \cite[p.\,187]{Skopets-Ponarin}. There are also several other useful properties of $4$-ball tetrahedra in \cite{Skopets-Ponarin}.

\begin{remark}
          \hspace{-0.15cm}{\bf .}
  \rm
By \cite[p.\,32]{Katsuura}, we have a very tight two-sided bound  (cf.~(\ref{eq:two-sided-estimate-path-tet}))
\begin{equation}
434.2^\circ \approx  2.41 \pi
\approx 1.5 \pi + 3 \, \arccos \frac{\sqrt{3}}{3} \leq \Sigma_T < 2.5 \pi = 450^\circ
\label{eq:two-sided-estim-cube-corner}
\end{equation}
for a general right-angled cube corner (trirectangular) tetrahedron $T$ with vertices $(0, 0, 0)$, $(a, 0, 0)$, $(0, b, 0)$, $(0, 0, c)$, where
$a$, $b$, $c$ are arbitrary positive numbers. Since three dihedral angles are equal to $90^\circ$, this rigidity follows.
The lower bound is attained if and only if $a=b=c$, i.e.~for a maximally symmetric cube corner tetrahedron
which is clearly a $4$-ball tetrahedron.
\label{rem:two-sided-estimates-cube-corner}
\end{remark}

\medskip

\noindent{\bf Example~1}
        \hspace{-0.15cm}{\bf .}
Another trivial example of a $4$-ball tetrahedron is the regular tetrahedron $T$. We observe that 
the geometrical mean of its inradius $r$ and circumradius $R$
is equal to the radius $\rho$ of its mid-sphere.
Indeed, if $a=1$ then $\vol_3 \, T=\sqrt{2}/12$ and by Remark~\ref{rm:karkas-radius}
we get
$$
rR=\frac{\sqrt6}{12}\cdot\frac{\sqrt6}{4}=\frac18=\Bigl(\frac{\sqrt2}{4}\Bigr)^2=\rho^2.
$$

\medskip

\noindent{\bf Example~2}
        \hspace{-0.15cm}{\bf .}
There exists also a $4$-ball tetrahedron none of whose faces is an equilateral triangle. For instance, if we prescribe the vertices
as follows
$$
A=(-2,0,0), B=(2,0,0), C=(0, \sqrt{21},0), D=(0, 5 \sqrt{21}/21, 4 \sqrt{105}/21),
$$
then the dihedral angles at egdes $|AC|=|BC|=5$   are $54.41^{\circ}$,
at $|AB|=|CD|=4$ are  $60.79^{\circ}$ and $83.62^{\circ}$, respectively,
and at $|AD|=|BD|=3$  are  $90^{\circ}$ exactly. Two faces $ACD$ and $BCD$ are Pythagorean right triangles and the other two are isosceles triangles.
We see that the sums of opposite edges are always $8$ and the sums of opposite angles are all equal to  $144.41^{\circ}$,
cf. properties 1) and  2) of Theorem~\ref{th:karkas-equivalent-definitions}. For this tetrahedron we have
$rR>\rho^2=9/5$. On the other hand, one can derive that $rR<\rho^2$ for a regular triangular pyramid
with height $h$, where $4h$ is the length of any edge of the triangular base.

According to \cite{Laszlo-2017}, we have
$$
R\ge \sqrt3\rho \ge 3r,
$$
where each of the two inequalities becomes an equality if and only if the $4$-ball tetrahedron is regular.

\medskip

In this section, let $a,b,c$ be the lengths of edges of some  triangular face of a given
tetrahedron~$T$. Let $d,e,f$ be the edges opposite to $a,b,c$, respectively.

\medskip

\noindent{\bf Example~3}
        \hspace{-0.15cm}{\bf .}
We can construct a $4$-ball tetrahedron with $a=2$, $d=20$ and
$b=c=e=f=11$. It has only two types
of its isosceles triangular faces:
$(2,11,11)$ with the largest angle $\arccos\frac{1}{11}\approx
84.78^\circ$ and
$(20,11,11)$ with the largest angle $2\arcsin\frac{10}{11}\approx
130.76^\circ$.
Hence, four of the corresponding six unfolded quadrilaterals from the tetrahedron
shell are non-convex, since
$84.78^\circ+130.76^\circ > 180^\circ$. Thus in general, there is no
inscribed
circle touching all four sides to the contrary to the situation 
sketched in Figure~\ref{fig:Kvant-Fig-7-8} (left).

\medskip

A tetrahedron cannot be constructed if one or two of its edges are long and the others are short. 
The edges $a$, $b$, $c$, $d$, $e$, $f$
have to satisfy all triangle inequalities for all four faces of $T$. However, this is
only a necessary condition as illustrated in the next example.

\medskip

\begin{figure}[htbp]
\centerline{
  \psfig{figure=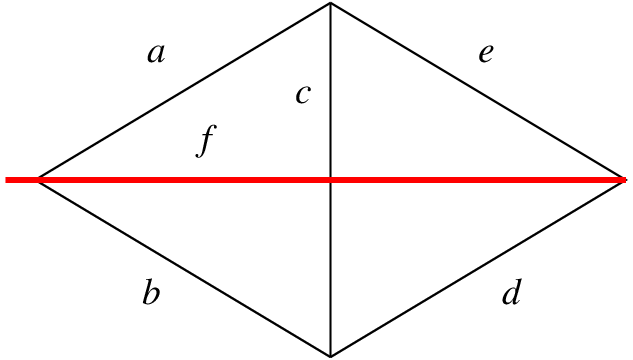, width=5.9cm,height=3.2cm}
  \qquad
  \qquad
    \psfig{figure=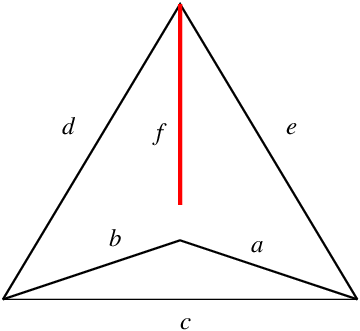, width=5.2cm,height=4.05cm}
}
\caption{Two typical cases when all triangle inequalities for all four faces
hold, but the Cayley-Menger determinant is negative.
Left: for $a=b=c=d=e=1$ the remaining thick red edge $f=7/4$ is too long
to construct a tetrahedron (see Example~4).
Right: for given edges $a=b=1.1$ and $c=d=e=2$ the remaining thick red
edge $f=1.1$
is too short to construct a ($4$-ball) tetrahedron, see Remark~\ref{rem:no-analog}.}
\label{fig:Cayley-Menger-negative}
\end{figure}

\noindent{\bf Example~4}
        \hspace{-0.15cm}{\bf .}
If $a=b=c=d=e=1$ and $f=7/4$, then all necessary triangle
inequalities are satisfied, but
the corresponding tetrahedron does not exist, because the altitude of
the adjacent equilateral triangular
faces $(a,b,c)$ and $(a,d,e)$ is $\sqrt{3}/2$, and $\sqrt{3}/2+\sqrt{3}/2<7/4$. So such a tetrahedron cannot
be constructed, see Figure~\ref{fig:Cayley-Menger-negative} (left) with more explanations why it is not possible.

\medskip

The next necessary and sufficient condition on six numbers to be the lengths of edges of some tetrahedron
is given in \cite[p.\,163]{Wirth-Dreiding}.

\begin{theorem}
          \hspace{-0.15cm}{\bf .}
The edges $a,b,c,d,e,f$ form a tetrahedron if and only if
\begin{itemize}

\item[i)] all edges of all four faces satisfy the corresponding triangle
inequalities and

\item[ii)] the Cayley-Menger determinant of the symmetric matrix
$$
  D_3=\det
\left(\begin{matrix} 0 & 1  & 1   & 1   & 1 \\
                        & 0  & a^2 & c^2 & e^2 \\
                        &    & 0   & b^2 & f^2 \\
                  &{\rm sym.}&     & 0   & d^2 \\
                        &    &     &     & 0
\end{matrix}\right)
$$
is positive. 
\label{th:sufficient-necessary-tet}
\end{itemize}
\end{theorem}

\begin{remark}
          \hspace{-0.15cm}{\bf .}
  \rm
  The Cayley-Menger determinant is closely related to the volume of $T$,
see e.g.   \cite[p.\,77]{Admal},
\cite[p.\,10]{BOOK}, \cite[p.\,434]{Fuhrmann},  
\cite[p.\,239]{Maehara-Martini}, \cite[p.\,737]{Menger} for various applications.
If $D_3=0$, then the associated
tetrahedron is degenerated with all vertices lying in one plane.
The conditions i) and ii) are independent. 
For the edges given in Example~4 we get a negative value of $D_3=-0.3828$.
We see that i) holds, but ii) does not hold. On the other hand, if $a=c=f=1$, $b=e=3$, and $d=5$,
then $a+c<b$ and $D_3=468$. Hence, i) does not hold, but ii) holds.

For the tetrahedra from Examples~1 and 2, we can derive that
$D_3=4$ for the regular tetrahedron with $a=1$
and
$D_3=10\,240$ for the tetrahedron with $a=b=5$, $c=f=4$, and $d=e=3$.
In Example~3, $D_3=256\,000$ and we get $D_3=8$ for the maximally symmetric
cube corner tetrahedron with $a=1$, cf.~Remark~\ref{rem:two-sided-estimates-cube-corner}.
\label{rm:CM-determinant-volume}
\end{remark}

\begin{remark}
          \hspace{-0.15cm}{\bf .}
  \rm
We show that the construction from Theorem~\ref{th:auxiliary-theorem} does not work in the case of tetrahedra.
Let $l_1, l_2, l_3, l_4$ be arbitrary positive numbers
and let
\begin{align*}
a=&l_1+l_2,\qquad b=l_1+l_3,\qquad c=l_2+l_3,\\
d=&l_3+l_4,\qquad e=l_2+l_4,\qquad f=l_1+l_4.
\end{align*}
Then
\begin{equation}
a+d=b+e=c+f,
\label{eq:no-working-for-tets}
\end{equation}
cf. condition 1) of Theorem~\ref{th:karkas-equivalent-definitions}.
We can prove like in the proof of Theorem~\ref{th:auxiliary-theorem}
that all necessary triangle inequalities for all triangular faces hold.
However, the
corresponding tetrahedron need not exist in general.
For instance, if $l_1=0.1$ and $l_2=l_3=l_4=1$, then
$a=b=f=1.1$, $c=d=e=2$ and
the Cayley-Menger determinant
is negative, $D_3=-11.84<0$, i.e., the condition ii) from Theorem~\ref{th:sufficient-necessary-tet}
is violated.
In other words, the tetrahedron with the above defined edges
cannot be constructed, see Figure~\ref{fig:Cayley-Menger-negative} (right).
This also means that for the three kissing balls with radii $l_1,l_2,l_3$ the fourth ball with radius $l_4$
can never kiss the first ball, i.e., the condition 4) from Theorem~\ref{th:karkas-equivalent-definitions} is violated.

Theorem~\ref{th:auxiliary-theorem} says that three balls can always be kissing,
which does not have to  be true for four balls. Theorem~\ref{th:karkas-equivalent-definitions} gives
necessary and sufficient conditions when this is possible.
\label{rem:no-analog}
\end{remark}

\begin{remark}
          \hspace{-0.15cm}{\bf .}
  \rm
According to \cite[p.\,424]{Blumenthal},
the Cayley-Menger determinant (in the notation of Theorem~\ref{th:sufficient-necessary-tet}) is equal to
\begin{align}
D_3=&\det
\left(\begin{matrix}  2a^2 & a^2+c^2-b^2 & a^2+e^2-f^2 \\
                            & 2c^2        & c^2+e^2-d^2 \\
                 {\rm sym.} &             & 2e^2
\end{matrix}\right)\nonumber\\[2mm]
=&2[4a^2c^2e^2+(a^2+c^2-b^2)(c^2+e^2-d^2)(a^2+e^2-f^2)\nonumber\\
  &-a^2(c^2+e^2-d^2)^2-c^2(a^2+e^2-f^2)^2-e^2(a^2+c^2-b^2)^2].     
\label{eq:D-3-Blumenthal}
\end{align}

\end{remark}

The proof of the next theorem is constructive.

\begin{theorem}
          \hspace{-0.15cm}{\bf .}
For any triangle there exist infinitely many 4-ball
tetrahedra having this triangle as one of its faces.
\label{th:karkas-from-face}
\end{theorem}

\begin{figure}[htbp]
\centerline{
  \psfig{figure=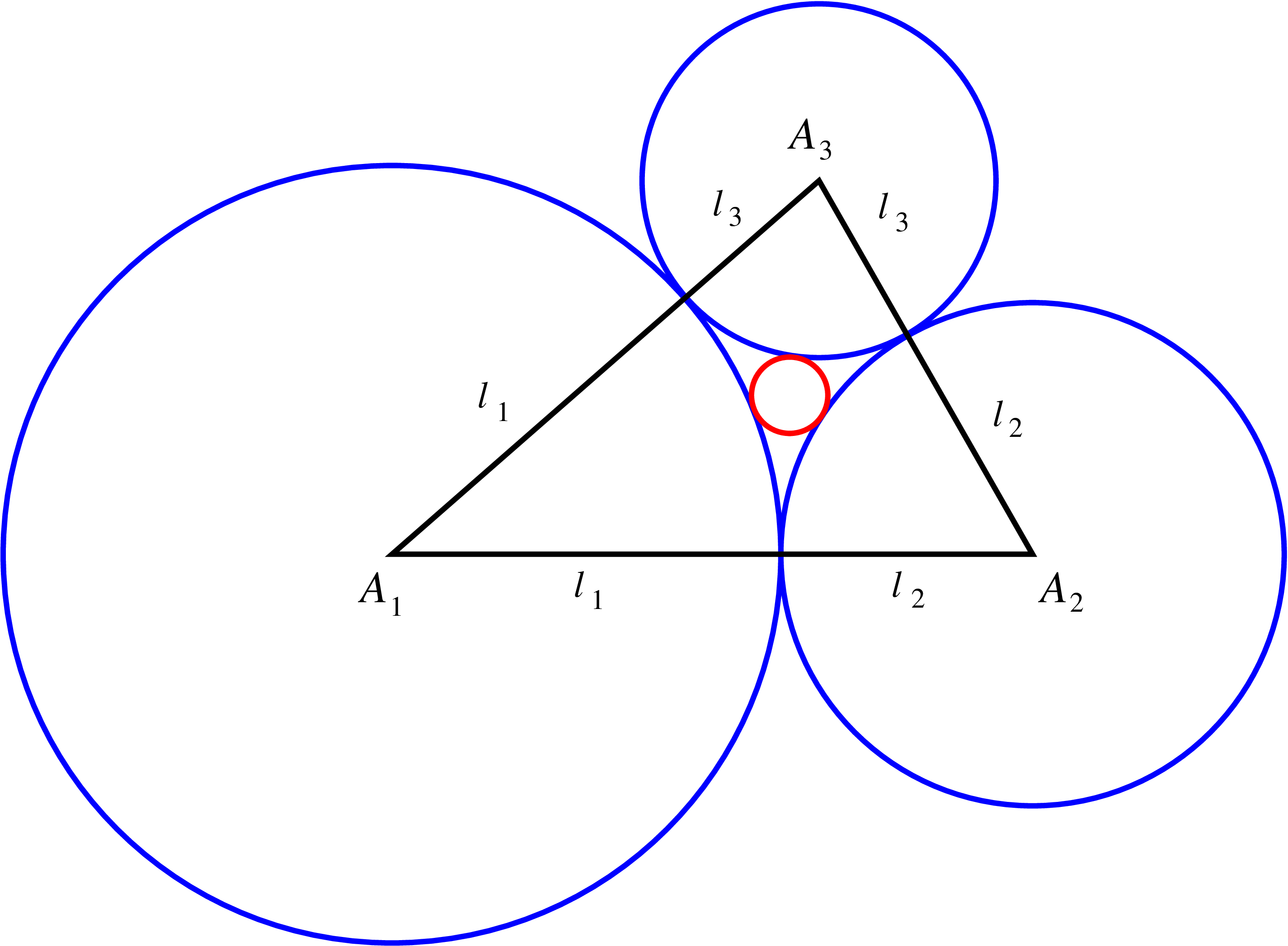, width=10.5cm,height=7.7cm}
}
\caption{
For the given three (blue) mutually externally touching circles
there exists a unique fourth (red) circle externally touching all three circles
(the 10th Apollonios problem).}
\label{fig:Soddy}
\end{figure}

\noindent{\bf P r o o f :}
Let $A_1A_2A_3$ be an arbitrary triangle in the plane $(x,y)$.
According to Theorem~\ref{th:auxiliary-theorem}, there exists a
unique correspondence between sides of $A_1A_2A_3$ and positive numbers
$l_1,l_2,l_3$. Hence, we can consider three circles with vertices $A_i$ and radii $l_i$, $i=1,2,3$,
introduced in Theorem~\ref{th:auxiliary-theorem}.
Then by Descartes' Theorem \cite{Lagarias} there there exists a uniquely
defined point $A_0$
inside $A_1A_2A_3$ which is the center
of the fourth circle (called the {\it Soddy circle} in \cite{Vandeghen})
which externally touches the three circles (see Figure~\ref{fig:Soddy}).
Hence, the interiors of all four circles are disjoint and each circle
touches the other three,
see Figure~\ref{fig:Soddy}. Denote by $l_0$ the radius of the Soddy circle (see \cite{Hirst} for its calculation).
According to Remarks~\ref{rm:CM-determinant-volume} and \ref{rem:no-analog}, we observe that $A_1A_2A_3A_0$ is
a degenerated $4$-ball tetrahedron.


Now consider three touching balls with the same centers at $A_1,A_2,A_3$
and same radii $l_1,l_2,l_3$.
Then for each sufficiently small $\varepsilon>0$ there exist exactly two
balls
with radius $l_4=l_0+\varepsilon$
touching the three given balls.
Then their centers $A_4$ and $\overline A_4$ are
symmetrically placed with respect to the plane $(x,y,0)$. Hence, by
property 4) of Theorem~\ref{th:karkas-equivalent-definitions},
$A_1A_2A_3A_4$ is a $4$-ball tetrahedron.
\halmos

\begin{remark}
          \hspace{-0.15cm}{\bf .}
          \rm
Let $a,b,c$ be sides of a fixed triangle and let us put
\begin{equation}
  d=k-a,\qquad e=k-b,\qquad f=k-c
  \label{eq:example-karkas-from-face}
\end{equation}
for some constant $k>\max(a,b,c)$ and substitute these values into (\ref{eq:D-3-Blumenthal}).
We again get relations (\ref{eq:no-working-for-tets}) and all triangle inequalities for all three remaining faces are satisfied, for instance,
$$
a+f > c-b + k-c = k-b = e.
$$
Moreover, we find by (\ref{eq:D-3-Blumenthal})
that $D_3=D_3(k)$ is a quadratic function in $k$ which is easy to investigate.

Setting
$$
c=1,\qquad X=a^2+b^2-1,\qquad Y=a^2+1-b^2,\qquad Z=b^2+1-a^2,
$$
we find by (\ref{eq:D-3-Blumenthal})   that
\begin{align*}
\tfrac12 D_3=&4a^2c^2e^2+(a^2+c^2-b^2)(c^2+e^2-d^2)(a^2+e^2-f^2)\\
     &-a^2(c^2+e^2-d^2)^2-c^2(a^2+e^2-f^2)^2-e^2(a^2+c^2-b^2)^2 \\[2mm]
=&4a^2(k-b)^2+Y(1+(k-b)^2-(k-a)^2)(a^2+(k-b)^2-(k-1)^2)\\
     &-a^2(1+(k-b)^2-(k-a)^2)^2-(a^2+(k-b)^2-(k-1)^2)^2-(k-b)^2Y^2
\\[2mm]
=&4a^2(k-b)^2+Y(b^2+1-a^2+2k(a-b))(a^2+b^2-1+2k(1-b))\\
     &-a^2(b^2+1-a^2+2k(a-b))^2-(a^2+b^2-1+2k(1-b))^2-(k-b)^2Y^2 \\[2mm]
=&4a^2(k^2-2bk +b^2)+Y(Z+2k(a-b))(X+2k(1-b))\\
     &-a^2(Z+2k(a-b))^2-(X+2k(1-b))^2-(k^2-2bk+b^2)Y^2 \\[2mm]
=&k^2[4a^2+4Y(a-b)(1-b)-4a^2(a-b)^2-4(1-b)^2-Y^2]\\
  &+k[-8a^2b+2XY(a-b)+2YZ(1-b)-4a^2Z(a-b)-4X(1-b)+2bY^2]\\
  &+4a^2b^2+XYZ-a^2Z^2-X^2-b^2Y^2 ,
\end{align*}
which is a quadratic function in $k$.
According to Theorem~\ref{th:karkas-from-face}, there exist infinitely
many $k$'s for which
the determinat $D_3=D_3(k)$ is positive for given $a$ and $b$.
This implies that the
corresponding discriminant has to be nonnegative. However, $D_3$ can be negative as shown in
Remark~\ref{rem:no-analog}.

The quadratic function $D_3=D_3(k)$ does not have to be positive
for all $k>\max(a,b,c)$.
To see this, consider an arbitrary fixed triangle $A_1A_2A_3$ and closed
balls $B_i$ centered
at $A_i$ with radii $l_i$. According to Theorem~\ref{th:auxiliary-theorem}, $l_1\leq l_2\leq l_3$
if and only if
$a\leq b\leq c$. Without loss of generality, we may assume that these
inequalities hold
and that $c=1$. If
$$
B_1\subset {\rm int conv}(B_2,B_3),
$$
then there is no ball $B_4$ for a sufficiently large radius $l_4$.
Hence, relations (\ref{eq:example-karkas-from-face})   cannot be satisfied.

\label{rem:D-3-quadratic}
\end{remark}

\medskip

In the next example, we illustrate that $4$-ball tetrahedra can be very badly shaped.

\medskip

\noindent{\bf Example~5}
        \hspace{-0.15cm}{\bf .}
We show that there exists a mirror symmetric $4$-ball tetrahedron
all of whose  four faces are obtuse triangles.
To satisfy condition 1) of Theorem~\ref{th:karkas-equivalent-definitions}    we put
$a=b=12$, $c=20$, $d=e=11$, and $f=3$.
Now from the Cosine Theorem we find that
all four faces are obtuse, since $12^2+12^2<20^2$, $11^2+11^2<20^2$, and
$3^2+11^2<12^2$ (twice).
All necessary triangle inequalities for all four faces are satisfied and by (\ref{eq:D-3-Blumenthal})
we have $D_3=448\,000 >0$.
So the assumptions of Theorem~\ref{th:sufficient-necessary-tet} hold.
Moreover, this $4$-ball tetrahedron has
three obtuse
dihedral angles $118.23^\circ$, $118.23^\circ$, $136.95^\circ$ at the
edges $d,e,f$,
respectively.

\medskip

From property 2) of Theorem~\ref{th:karkas-equivalent-definitions}
we see that no path tetrahedron is $4$-ball,
since  the two dihedral angles at
opposite edges $AC$ and $BO$ in Figure~\ref{fig:path-tet} are right and the sum
$\Sigma_T$ of all six dihedral angles is less than $3 \pi$ by (\ref{eq:two-sided-estimate-general-tet}). 
Anyway, we have the following theorem.





\begin{theorem}
          \hspace{-0.15cm}{\bf .}
  Let $T$ be a $4$-ball tetrahedron. Let the center of the mid-sphere lie in the interior of $T$. Then $T$ can be partitioned face-to-face into $24$
  path tetrahedra.
  \label{th:karkas-into-24-path-tets}
\end{theorem}

\begin{figure}[htbp]
\centerline{
  \psfig{figure=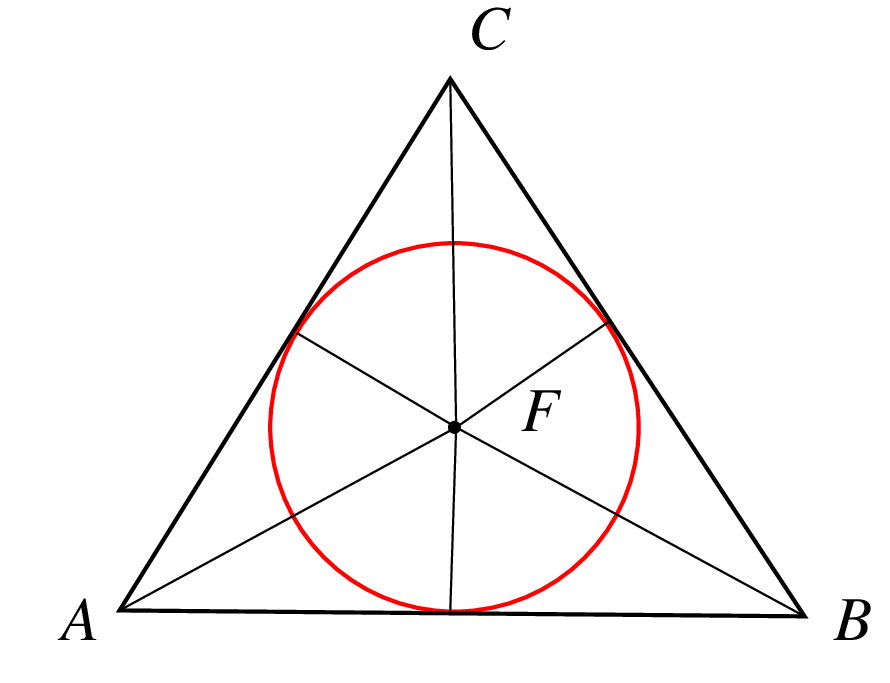, width=5.5cm,height=4.15cm}
  \qquad
    \psfig{figure=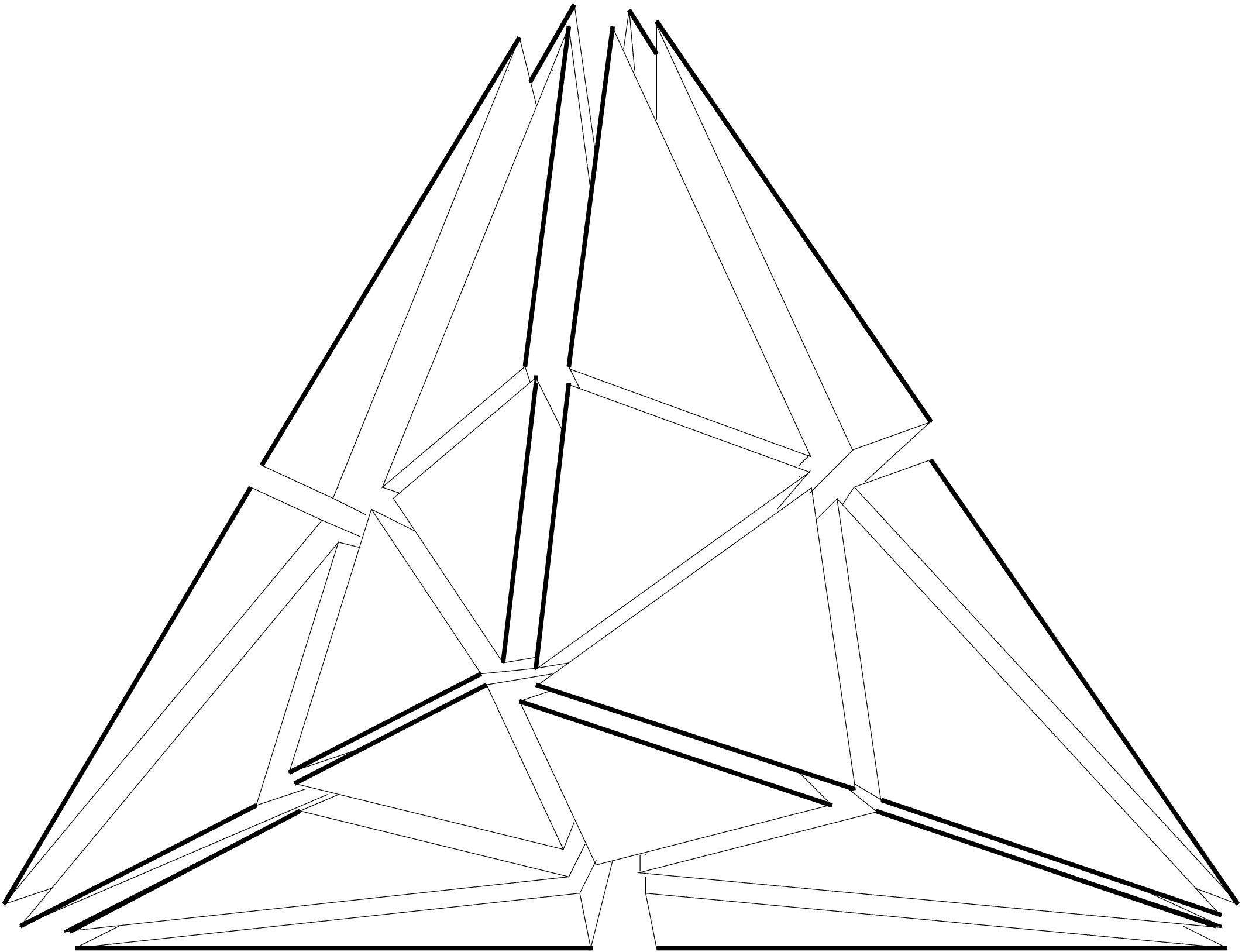, width=8cm,height=5.4cm}
}
\caption{
Left: Partition of the face $ABC$ in $6$ right triangles using the incenter $F$. Right: Face-to-face partition of the $4$-ball
tetrahedron from Theorem~\ref{th:karkas-into-24-path-tets} into $24$ path tetrahedra. In this way all obtuse dihedral angles
disappear due to dissections. Moreover, every path tetrahedron can be further
partitioned face-to-face into $8$ smaller path subtetrahedra by a red refinement algorithm, see Figure~\ref{fig:path-tet} (right)
and \cite[p.\,90]{BOOK}.}
\label{fig:24-path-tets}
\end{figure}

\noindent{\bf P r o o f :}
Let $S$ be the  mid-sphere of a given $4$-ball tetrahedron $ABCD$. Each edge can be thus divided into two straight line segments
$l_i$ and $l_j$, $i\neq j$, $i,j\in\{1,2,3,4\}$, by the corresponding touching points (see Remark~\ref{rm:karkas-radius}).
The face $ABC$ can be thus decomposed into six right triangles 
whose common vertex $F$ is the center of the inscribed circle to $ABC$, 
see Figure~\ref{fig:24-path-tets} (left).
The inscribed circle is in fact the intersection
$S\cap ABC$. Denoting by $G$ the center of the mid-sphere $S$, it is clear that $FG$ is perpendicular to $ABC$,
see property 5) of Theorem~\ref{th:karkas-equivalent-definitions}.

By taking the convex hulls of the six right triangles with center $G$ we define six path tetrahedra.  Since edges of adjacent faces are divided
in the same way,  see Figures~\ref{fig:karkas-tangents} and \ref{fig:Kvant-Fig-7-8} (right),
we can repeat the above construction also for the remaining three faces $ABD$, $ACD$ and $BCD$.
  In this way we get a face-to-face partition of the $4$-ball tetrahedron $ABCD$ into 24 path subtetrahedra,
  see Figure~\ref{fig:24-path-tets} (right). \halmos

\begin{remark}
          \hspace{-0.15cm}{\bf .}
  \rm
If $G$ from the above proof lies on the boundary $\partial T$,
then the $4$-ball tetrahedron $T$ can be partitoned face-to-face into $18$
path tetrahedra, since the six remaining path tetrahedra (out of $24$) are
degenerated. Obviously, the center $G$ of $S$ cannot lie on edges of
$T$, and thus we cannot have less than $18$ tetrahedra in the construction.
If $G\in T$, then during the refinement process into
path tetrahedra all
possible obtuse dihedral angles disappear, since path tetrahedra are
always nonobtuse. These special tetrahedra play an important role in
numerical
analysis, since they produce monotone stiffness matrices (see \cite[p.\,118]{BOOK}
for the discrete maximum principle).
\end{remark}

\begin{remark}
          \hspace{-0.15cm}{\bf .}
  \rm
In \cite[p.\,99]{BOOK} we present a similar algorithm as in the
proof of Theorem~\ref{th:karkas-into-24-path-tets} that yields the so-called face-to-face yellow refinement in
which  all path subtetrahedra contain the center of the circumscribed sphere (cf. Example~6 below). Then there are no
obtuse angles which enables us to satisfy the discrete maximum principle of several problems of mathematical
physics solved by the finite element method, see e.g. \cite[p.\,115]{BOOK} and \cite{Vejchodsky}.
\end{remark}

\section{Two-sided estimation of the dihedral angle sum for $4$-ball tetrahedra}
\label{sec:estimation-karkas}

First we show that three dihedral angles associated with  only one vertex uniquely determine the entire $4$-ball tetrahedron
(up to scaling). Consider three half-lines starting
at one point that do not lie in one plane. Their convex hull will be called a
{\it triangular cone}. The proof of the next theorem is again constructive.


\begin{theorem}
          \hspace{-0.15cm}{\bf .}
  For any triangular cone there exists exactly one $4$-ball tetrahedron up to scaling such that its
  three adjacent edges are contained in the edges of the cone. 
\label{th:cone-karkas}
\end{theorem}

\noindent{\bf P r o o f :}
Let $A_1=(0,0,0)$ be the vertex of a given triangular cone. Its three faces
    are contained in the planes
    \begin{equation}
    u_jx+v_jy+w_jz=0, \qquad j=2,3,4,                     \label{eq:plane}
    \end{equation}
 where $n_j=(u_j,v_j,w_j)$ are the associated outward unit normals. 
 Any couple of these equations determines one dihedral angle $\arccos(-n_i\cdot n_j)$, $i\neq j$, at one of the three edges of the cone,
 i.e. three half-lines starting at $A_1$. Let 
    $$
    (x-x_0)^2+(y-y_0)^2+(z-z_0)^2=1
    $$    
    be the unit sphere touching the three non-coplanar edges, where $G=(x_0,y_0,z_0)$ is its center.
    The intersection of this sphere with planes (\ref{eq:plane}) uniquely determines
    there circles as marked in Figure~\ref{fig:karkas-tangents}
    The three touching points are indicated by three bullets on edges emanating from $A_1$. 
    
    There exists exactly one plane which touches all three circles at another three points
    (see the remaining three bullets in Figure~\ref{fig:karkas-tangents})
    so that the circles
    belong to the half-space containing $A_1$. This plane uniquely determines
    the points $A_2, A_3$ and $A_4$ on the corresponding half-lines. In this way we obtain
    the $4$-ball tetrahedron $A_1A_2A_3A_4$, since the sums of lengths of opposite edges
    are all equal to $l_1+l_2+l_3+l_4$. It is obvious that the above unit sphere is the corresponding mid-sphere.
    For a nonunit sphere suitable scaling arguments can be used.
\halmos

\medskip

  The opposite statement is obvious, since any given $4$-ball tetrahedron can be inserted to
     an infinite triangular cone whose three edges contain three adjacent edges of the tetrahedron.

\medskip

\noindent{\bf Example~6}
        \hspace{-0.15cm}{\bf .}
    All $4$-ball tetrahedra from Examples~1, 2, 3 and 5 are mirror symmetric.
    Consider now an asymmetric $4$-ball tetrahedron $T$ all of whose edges have different lengths:
    \begin{align*}
    a=|BC|=5,~\qquad &b=|AC|=12,\qquad c=|AB|=13,\\
    d=|AD|=11,\qquad &e=|BD|=4,~\qquad f=|CD|=3,
    \end{align*}
    where $A=(12,0,0)$, $B=(0,5,0)$, $C=(0,0,0)$, $D=(\tfrac43,\tfrac95,z)$  are vertices of $T$  
    and $z=8\sqrt{14}/15$. Note that  its faces $ABC$ and $BCD$ are Pythagorean right triangles
    and vol$_3(T)=10z$. We further note that $T$ does not contain its circumcenter $(6,5/2,-15/\sqrt{14})$, since
    its last coordinate is negative. This means that we cannot apply Theorem~8.8 from \cite[p.\,98]{BOOK} for partition
    of $T$ into $24$ path tetrahedra. Note that $T$ has 3 obtuse dihedral angles at edges containing the vertex $D$.   
    By Theorem~\ref{th:auxiliary-theorem}
    we find that $l_1=2$, $l_2=10$, $l_3=3$, and thus $l_4=1$.
    Moreover, by Remark~\ref{rm:karkas-radius} 
    the radius of the mid-sphere is $\rho=4/z$. Its center $G=(2,2,1/\sqrt{56})$ lies inside $T$, and thus we can apply
    the above Theorem~\ref{th:karkas-into-24-path-tets} to decompose $T$ face-to-face into $24$ path tetrahedra.

    Taking $A_1=C$ in the proof of     Theorem~\ref{th:cone-karkas}, 
    we find that circles inscribed to the faces $ABC$, $ACD$ and $BCD$
    have radii $2$, $\sqrt{260}/13$ and $1$, respectively. The fourth face $ABD$ touches
    these three circles as well. 

\medskip

The upper bound of the  dihedral angle sum for $4$-ball tetrahedra is larger than for path tetrahedra (see Theorem~\ref{th:sum-for-path-tets}). For instance,
$\Sigma_T \approx 479.58^\circ \approx 2.66\pi > 2.5 \pi$ for the $4$-ball tetrahedron from Example~3
and $\Sigma_T=453.33^\circ$ for the asymmetric $4$-ball tetrahedron from Example~6.  Actually, we have the following
two-sided bound, where     the lower bound is approximately $423.17^\circ\approx 2.35\pi$.

\begin{theorem}
          \hspace{-0.15cm}{\bf .}
Let $T$ be a $4$-ball tetrahedron. Then
\begin{equation}
  6 \, \arccos\frac13 \leq \Sigma_T <3\pi
\label{eq:tw-sided-estim-karkas}
\end{equation}
and each value in this interval can be attained by some $4$-ball tetrahedron.
\label{th:main-estimate-karkas-tet}
\end{theorem}

\noindent{\bf P r o o f :}
First we show that the upper bound $3\pi$ is tight. Consider
a one-parameter family of tetrahedra with an equilateral triangular
base and height $h > 0$ as a parameter. Let the remaining vertex (apex) orthogonally
projects to the center of gravity $F$ of the base, see Figure~\ref{fig:24-path-tets} (left).
Moving the apex to
$F$, the three dihedral angles at the base edges converge to zero and the
remaining three dihedral angles converge to $\pi$ as $h\to 0$. In this case, the centers of mid-spheres lie outside tetrahedra.

We can easily check that the lower bound (\ref{eq:tw-sided-estim-karkas}) is then reachable for the regular
tetrahedron all of whose dihedral angles are arccos$\tfrac13\approx 70.5^\circ$.
The function $h\mapsto \Sigma_T (h)$ is decreasing on the interval      $(0,h_0]$ and increasing on $[h_0,\infty)$, where $h_0$ is the
    height of the regular tetrahedron. Moreover, each value from interval (\ref{eq:tw-sided-estim-karkas}) is reachable,
since $\Sigma_T$ changes continuously with $h$  in the above one-parameter family.

Finally, we show that there is no $4$-ball tetrahedron $T$
with $\Sigma_T$ smaller than $6 \, \arccos\tfrac13$. By \cite{Gad}, the dihedral angle and the angle between  outward normals
to the faces add up to $\pi$. Thus, the angle between outward unit normals to two faces  is uniquely associated with
the geodesic distance along a great circle on the unit sphere. According to
property 2) of Theorem~\ref{th:karkas-equivalent-definitions}, the sums of dihedral angles to opposite edges
are equal. Therefore, the sums of geodesic distances corresponding to
opposite edges are also equal.
Consequently, we shall consider a spherical triangle
with sides $a,b,c\in(0,2\pi)$ and an unknown angle $\alpha$ opposite to $a$. Then
the Cosine Theorem for the unit sphere (i.e. the spherical law of cosines)
reads, see e.g. \cite[p.~85]{Rektorys},
\begin{equation}
\cos a=\cos b\cos c+\sin b\sin c\cos\alpha.                 \label{10}
\end{equation}
Note that the symbols $a,b,c$
have now different
meanings than those in Section~\ref{sec:basic-properties-karkas}.

\begin{figure}[htbp]
\centerline{
  \psfig{figure=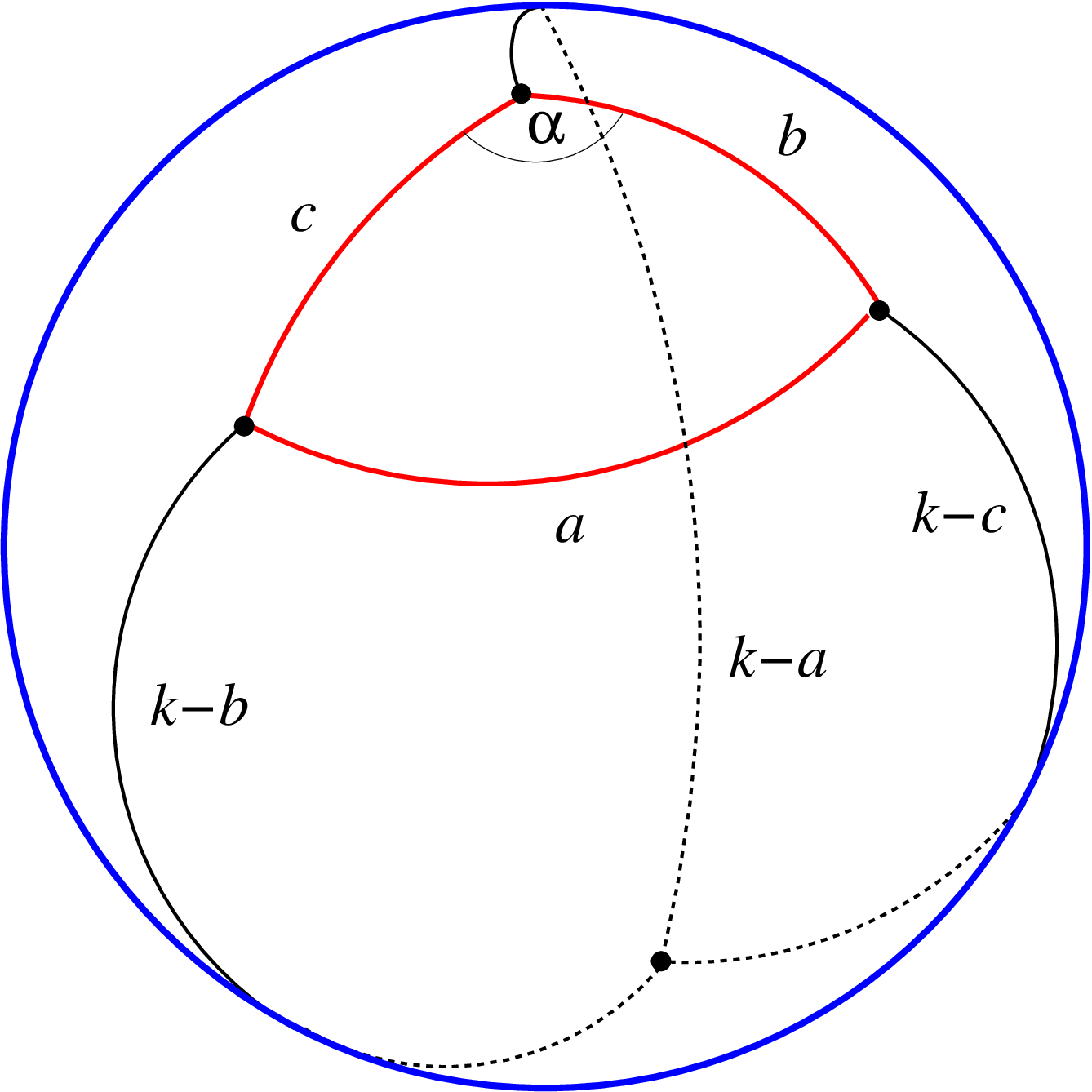, width=6cm,height=6cm}
}
\caption{Schematic illustration of the geodesic distances $a,b,c,k-a,k-b,k-c$ between four outward normals (marked by
bullets) to the faces of a 4-ball tetrahedron on the unit sphere.}
\label{fig:sphere-distances}
\end{figure}

If three dihedral angles and their three corresponding geodesic distances $a,b,c\in(0,2\pi)$ are given, then by Theorem~\ref{th:cone-karkas}
there exists a uniquely determined $4$-ball tetrahedron $T$ (up to scaling) such that $k-a$, $k-b$, $k-c$ are opposite sides to $a,b,c$, respectively.
Here $k>\max(a,b,c)$ is an unknown constant which is also unique for fixed $a, b, c$. According to \cite{Gad},
\begin{equation}
  \Sigma_T+\Gamma_T=6\pi,
\label{11}
\end{equation}
where $\Gamma_T$ is the sum of all six geodesic distances between four outward normals marked by
four  bullets in Figure~\ref{fig:sphere-distances}. According to property 2) of Theorem~\ref{th:karkas-equivalent-definitions},
it is obvious that for a $4$-ball tetrahedron we have

\begin{equation}
\Gamma_T=a+b+c+(k-a)+(k-b)+(k-c)=3k.
\label{eq:Gamma_T}
\end{equation}

\bigskip

To introduce the main idea of the proof, we will first assume for simplicity that $b=c$, i.e., we shall investigate an isosceles
spherical triangle.  Then (\ref{10}) reduces to
$$
\cos a=\cos^2b+\sin^2b\cos\alpha
$$
and we have
\begin{equation}
\alpha=\arccos\frac{\cos a-\cos^2 b}{\sin^2 b}.                   \label{12}
\end{equation}
Now consider two neighboring spherical triangles with sides $b$, $k-a$ and $k-b$. 
The angle between the curved sides $b$ and $k-a$ is equal to $\pi-\tfrac12\alpha$
due to symmetry, since $\alpha+(\pi-\tfrac12\pi)+(\pi-\tfrac12\pi)=2\pi$.
Then, similarily to (\ref{10}), we get by the Cosine Theorem for sides
\begin{align*}
\cos(k-b)&=\cos b\cos(k-a)+\sin b\sin(k-a)\cos(\pi-\tfrac12\alpha)\\
        &=\cos b\cos(k-a)-\sin b\sin(k-a)\cos\tfrac{\alpha}{2}.
\end{align*}
Applying (\ref{12}) and simplifying by setting
$$
C=\cos a,\quad S=\sin a, \quad B=\cos b, \quad D= \sin b, 
\quad A=\cos\Bigl(\frac12\arccos\frac{C-B^2}{D^2}\Bigr),
$$
we get
$$
B\cos k+D\sin k=B(C\cos k+S\sin k)-AD(C\sin k - S\cos k).
$$
Hence,
$$
(B-BC-ADS)\cos k = (BS-D-ACD)\sin k
$$
and thus 
we can express the unknown value $k$ explicitly using the given values $a$ and $b=c$,
namely,
$$
k = \arctan \frac{B-BC-ADS}{BS-D-ACD}.
$$
By differentiating with respect to $a$ and $b$, we find that $k = k(a,b)$ attains its global maximum value over $(0,2\pi)\times(0,2\pi)$ when 
$a=b=\pi-\arccos\tfrac13$. In this case,
$\alpha=\tfrac23 \pi =120^\circ$,
\begin{equation}
k_{\rm max}=2\,\arccos\Bigl(-\frac13\Bigr)=2a=3.821266\dots (\approx 218.942441^\circ) \label{13}
\end{equation}
and
$$
\quad A=\frac12,\quad B=C=-\frac13,\quad D=S=-\frac{2\sqrt{2}}{3},
\quad \cos k=-\frac79,\quad \sin k=-\frac{4\sqrt{2}}{9}.
$$
According to (\ref{11}), maximizing $\Gamma_T$ is obviously equivalent to minimizing $\Sigma_T$.

\medskip

A tedious calculation can be done similarly for the case $b\neq c$ which
leads to
the same maximum value (\ref{13}) of $k=k(a,b,c)$ over the cube $(0,2\pi)\times(0,2\pi)\times(0,2\pi)$. 
The lower bound in (\ref{eq:tw-sided-estim-karkas}) now follows from (\ref{11}) and (\ref{eq:Gamma_T}).     \halmos

\begin{remark}
          \hspace{-0.15cm}{\bf .}
  \rm
  An analogous two-sided bound for the dihedral angle sum of equifacial tetrahedra all of whose faces are congruent triangles
  (see \cite[p.~33]{Katsuura}) is
\begin{equation}  
 2 \pi < \Sigma_T  \le 6 \, \arccos  \frac{1}{3}.  \label{14}
\end{equation}
Hence, this estimate is complementary to estimate (\ref{eq:tw-sided-estim-karkas}). The intersection of estimates
(\ref{eq:tw-sided-estim-karkas}) and
(\ref{14}) corresponds to the regular tetrahedron only. 
\end{remark}

\bigskip

\bf Acknowledgments. \rm The authors are indebted to Jan Chleboun,  Leandro Farina, Ivan Gudoshnikov, Lawrence Somer, and
Tom\'a\v s Vejchodsk\'y for fruitful discussions.
Supported by the Czech Science Foundation No. 24-10586S  and the Czech Academy of Sciences (RVO 67985840).

\renewcommand\refname

\end{document}